\documentclass[12pt]{amsart} 
\usepackage{enumerate}
\input{amssymb.sty}
\title[Non-amenability of $\mathcal{B}(\ell_p)$]
{A note on non-amenability of $\mathcal{B}(\ell_p)$ for $p=1,2$} 
\author{Narutaka Ozawa} 
\address{Department of Mathematical Science, 
University of Tokyo, 153-8914, Japan } 
\email{narutaka@ms.u-tokyo.ac.jp} 

\newtheorem{thm}{Theorem}[section]
\newtheorem{cor}[thm]{Corollary} 
\newtheorem{lem}[thm]{Lemma} 
\theoremstyle{definition} 
\newtheorem{rem}[thm]{Remark}
\newtheorem{defn}[thm]{Definition} 
\newcommand{\N}{{\mathbb N}} 
\newcommand{\Z}{{\mathbb Z}} 
\newcommand{\R}{{\mathbb R}} 
\newcommand{\C}{{\mathbb C}} 
\newcommand{\F}{{\mathbb F}} 
\newcommand{\B}{\mathcal{B}}
\newcommand{\hh}{\mathcal{H}}
\newcommand{\p}{\varphi} 
\newcommand{\G}{\Gamma} 
\newcommand{\Ll}{\Lambda_l} 
\newcommand{\Tr}{\mathop{\mathrm{Tr}}} 
\newcommand{\Lip}{\mathop{\mathrm{Lip}}} 
\newcommand{\reg}{\mathrm{reg}} 
\newcommand{\supp}{\mathop{\mathrm{supp}}} 
\newcommand{\Sg}{\mathcal{S}} 
\newcommand{\e}{\varepsilon} 
\newcommand{\dist}{\mathop{\mathrm{dist}}} 
\begin{document}
\begin{abstract}
This is an expository note on non-amenabilty of 
the Banach algebra $\B(\ell_p)$ for $p=1,2$. 
These were proved respectively by Connes and Read 
via very different methods. 
We give a single proof which reproves both.
\end{abstract}
\maketitle
\section{Introduction}
Johnson \cite{johnson} introduced the notion 
of amenability for Banach algebra 
and asked whether $\B(X)$, the Banach algebra of 
bounded linear map on a Banach space $X$, 
can be amenable. 
This problem is still widely open. 
For the special case where $X=\ell_p$, 
it was proved by Connes \cite{connes} that $\B(\ell_2)$ 
and $\prod_{n=1}^\infty \B(\ell_2^n)$ are not amenable, 
and more recently, by Read \cite{read} that 
$\B(\ell_1)$ and $\prod_{n=1}^\infty \B(\ell_p^n)$ 
for $p\in [1,\infty]\setminus\{2\}$ are not amenable. 
It is left open whether $\B(\ell_p)$ is amenable or not 
for $p\in(1,\infty)\setminus\{2\}$. 
A somewhat simplified proof of Read's results 
was given by Pisier \cite{pisier} where it was observed 
that the proof also gives 
the non-amenability of `regular' $\B(\ell_p)$. 
In this note, we give a single proof 
which reproves all above results. 
The proof is inspired from Read's approach 
and was suggested by Pisier. 
\begin{thm}\label{main}
The following Banach algebras are not amenable; 
$\B(\ell_p)$ for $p\in\{1,2,\infty\}$, 
$\prod_{n=1}^\infty \B(\ell_p^n)$ for $p\in[1,\infty]$ and 
$\B_r(\ell_p)$ for $p\in[1,\infty]$. 
\end{thm}
Here, $\B_r(\ell_p)$ is the Banach algebra 
consisting of those operators $x\in \B(\ell_p)$ with 
\[
\| x\|_{\B_r(\ell_p)} := \|\, |x|\, \|_{\B(\ell_p)}<\infty, 
\]
where $|x|=(|x_{ij}|)_{i,j}$ for $x=(x_{ij})\in\B(\ell_p)$. 
We use the following definition of amenability 
of Banach algebras, which is equivalent to the standard one 
by Johnson's theorem \cite{johnson2}. 
A general reference for amenable Banach algebra is \cite{runde}. 
\begin{defn}\label{def}
A unital Banach algebra $A$ is \emph{amenable} 
if there exists a constant $C>0$ with the following property; 
for any finite set $F\subset A$ and $\e>0$, 
there exists $T=\sum_{i=1}^r a_i\otimes b_i\in A\otimes A$ 
such that 
\begin{enumerate}[(i)]
\item 
$\mathrm{prod}(T)=\sum_{i=1}^r a_ib_i=1,$
\item
$\|x\cdot T-T\cdot x\|_{\wedge}
=\| \sum_{i=1}^r xa_i\otimes b_i 
- \sum_{i=1}^r a_i\otimes b_ix \|_{\wedge}\le\e$ 
for $x\in F$, 
\item
$\| T\|_{\wedge}\le C$, 
\end{enumerate}
where $\|\,\cdot\,\|_{\wedge}$ denotes the projective tensor 
norm on $A\otimes A$. 
\end{defn}

It turns out that the condition (iii) 
is irrelevant to Theorem \ref{main}. 

\section{Some Facts on $\ell_p$}
In this section, we prove some elementary facts on $\ell_p$. 
Fix $p\in[1,\infty]$ and let $q$ be so that $1/p+1/q=1$. 
We denote by $\{\delta_m\}$ (resp.\ $\{\delta_m^*\}$) 
the standard basis of $\ell_p$ (resp.\ $\ell_{q}$). 
We temporally denote by $\|x\|_{\reg}=\|\,|x|\,\|$ 
for $x$ in $\B(\ell_p)$ or $\B(\ell_q)$. 

\begin{lem}\label{lpineq}
For $x\in \B(\ell_p,\ell_p^N)$ and $y\in \B(\ell_q,\ell_q^N)$, 
we have 
\begin{enumerate}[(i)]
\item 
$\sum_{m=1}^\infty \|x\delta_m\|_q\|y\delta^*_m\|_p\le N\|x\|\,\|y\|$, 
\item
$\sum_{m=1}^\infty \|x\delta_m\|_p\|y\delta^*_m\|_q
\le N\|x\|_{\reg}\|y\|_{\reg}$, 
\item
$\sum_{m=1}^\infty \| x\delta_m\|_2\|y\delta^*_m\|_2
\le N(\|x\|\|x\|_{\reg}\|y\|\|y\|_{\reg})^{1/2}$. 
\end{enumerate}
\end{lem}
\begin{rem}
The estimate 
$\sum_{m=1}^\infty \| x\delta_m\|_2\|y\delta^*_m\|_2
\leq N\|x\|\|y\|$ needs not be true (even up to constant multiple) 
for $x\in \B(\ell_p,\ell_p^N)$ and $y\in \B(\ell_q,\ell_q^N)$. 
However, it is not hard to see (from the following proof) 
this is true if $p\in\{1,2\}$ or if 
$x\in \B(\ell_p^N)$ and $y\in \B(\ell_q^N)$. 
\end{rem}
\begin{proof}
The part (i) follows from H\"{o}lder inequality; 
\begin{align*}
\sum_{m=1}^\infty \|x\delta_m\|_q\|y\delta^*_m\|_p
&\le (\sum_{m=1}^\infty \|x\delta_m\|_q^q)^{1/q}
 (\sum_{m=1}^\infty \|y\delta^*_m\|_p^p)^{1/p}\\
&= (\sum_{k=1}^N \|x^*\delta^*_k\|_q^q)^{1/q}
 (\sum_{k=1}^N \|y^*\delta_m\|_p^p)^{1/p}\\
&\le N\|x\|\,\|y\|.
\end{align*}
For the part (ii), 
we may assume $x=|x|$. 
Then, by the duality $\ell_q=(\ell_p)^*$, we have 
\begin{align*}
(\sum_{m=1}^\infty \|x\delta_m\|_p^q)^{1/q}
&= \sup_\lambda \sum_{m=1}^\infty \lambda_m\|x\delta_m\|_p 
\le \sup_\lambda \sum_{m=1}^\infty \lambda_m\|x\delta_m\|_1\\ 
&= \sup_\lambda \| x\lambda\|_1 
\le N^{1/q}\sup_\lambda \|x\lambda\|_p\\
&\le N^{1/q}\| x\|_{\reg}, 
\end{align*}
where the supremum is taken over all $\lambda=(\lambda_m)_m$ 
with $\|\lambda\|_p\le 1$. 
The rest of the proof of part (ii) is same as (i). 
Finally, the part (iii) follows from (i) and (ii) 
combined with the Cauchy-Schwarz inequality and the fact 
$\|\lambda\|_2\le(\|\lambda\|_p\|\lambda\|_q)^{1/2}$. 
\end{proof}
\section{Some Facts on $\mathrm{SL}(3,\Z)$}\label{sl}
In this section, we collect some known 
facts on the group $\mathrm{SL}(3,\Z)$. 
For notational simplicity, we denote $\mathrm{SL}(3,\Z)$ by $\G$. 
The group $\G$ has Kazhdan's property $\mathrm{(T)}$. 
In particular, it is finitely generated. 
(In fact, it is generated by two elements.) 
A general reference for Kazhdan's property $\mathrm{(T)}$ is \cite{hv}. 
\begin{thm}\label{prt} 
The group $\G=\mathrm{SL}(3,\Z)$ has 
Kazhdan's property $\mathrm{(T)}$. 
Namely, for any finite set of generators $\Sg$ of $\G$, 
there exists a positive number $R=R(\G,\Sg)$ 
with the following property; 
for any unitary representation $\G$ on a Hilbert space $\hh$ 
and a vector $\xi\in\hh$, 
there exists a $\pi$-invariant vector $\eta\in\hh$ such that 
\[
\|\xi-\eta\|\le R\max_{g\in \Sg}\|\pi(g)\xi-\xi\|.
\]
\end{thm}
This theorem is non-trivial only when 
$\max_{g\in \Sg}\|\pi(g)\xi-\xi\|<\|\xi\|/R$ 
since the zero vector is $\pi$-invariant. 
See Appendix for a related fact on 
actions on Banach spaces.

We next present explicit examples of actions of $\G$. 
For each prime number $l$, let $\F_l=\Z/l\Z$ be the finite field 
with $l$ elements and define the ``finite projective plane'' 
$\Ll$ by $\Ll=(\F_l^3-\{0\})/\F_l^\times$. 
The action $\G$ on $\F_l^3$ 
(through linear transformation by $\mathrm{SL}(3,\F_l)$) 
induces an action $\sigma_l$ on the set $\Ll$. 
It is not hard to see that $\sigma_l$ is 2-transitive, 
i.e., for any $s_{ij}\in \Ll$ ($i,j\in\{1,2\}$) 
with $s_{1j}\neq s_{2j}$ for $j\in\{1,2\}$, 
there exists $g\in\G$ such that 
$\sigma_l(g)(s_{i1})=s_{i2}$ for $i=1,2$. 
In other words, the product action 
$\sigma_l\times\sigma_l$ of $\G$ on $\Ll\times\Ll$ 
has exactly two orbits (diagonals and off-diagonals). 
For $p\in[1,\infty]$, we denote by $\pi_l^p$ the 
isometric action of $\G$ on $\ell_p(\Ll)$ induced 
from $\sigma_l$. 
Let $\Lambda=\bigsqcup_{l\mbox{ prime}}\Ll$ be 
the disjoint union of $\Lambda_l$'s and 
let $\pi^p=\bigoplus\pi_l^p$ be the isometric action 
of $\G$ on $\ell_p(\Lambda)$. 
For a later use, we choose a subset $C_l\subset\Ll$ 
with $|C_l|=(l^2+l)/2=(|\Ll|-1)/2$ 
and let $v^p_l\in \B(\ell_p(\Ll))$ be the 
isometry given by $\upsilon^p_l\delta_s=\delta_s$ if $s\in C_l$ 
and $\upsilon^p_l\delta_s=-\delta_s$ if $s\notin C_l$. 
Finally, set 
$\upsilon^p=\bigoplus\upsilon^p_l\in\B(\ell_p(\Lambda))$. 
\section{Some Facts on $\ell_2\hat{\otimes}\ell_2$}\label{mazur}
In this section, we review a well-known fact
about operator inequality. 
A general reference is \cite{bhatia}. 
Let $S_1$ (resp.\ $S_2$) be the 
Banach space of trace (resp.\ Hilbert-Schmidt) 
class operators on $\ell_2$. 
The `noncommutative' Mazur map $\p\colon S_1\to S_2$ is defined 
by $\p(T)=U|T|^{1/2}$, where $T=U|T|$ is the polar decomposition 
and modulus and square root are taken 
as an operator on a Hilbert space. 
By the uniqueness of polar decomposition, 
we have $\p(USU^*)=U\p(S)U^*$ for any unitary operator $U$.
As in the case of ordinary Mazur map, 
we have the following.

\begin{thm}
The noncommutative Mazur map $\p\colon S_1\to S_2$ 
is a uniform homeomorphism between the unit ball of 
$S_1$ and that of $S_2$. 
\end{thm}
\begin{proof}
We give a proof for the reader's convenience. 
Temporally, we denote by $\|\,\cdot\,\|_p$ 
the Schatten $p$-norm. 
It suffices to show that $\p$ is a uniform homeomorphism 
between the spheres. 
Let $S=U|S|$ and $T=V|T|$ be given so that 
$\|S\|_1=1=\|T\|_1$. Let $\e=\|S-T\|_1$. 
Then, $\| |S|-|T|\|_1\le2\e^{1/2}$ by \cite{kosaki} 
and $\| |S|^{1/2}-|T|^{1/2}\|_2\le2\e^{1/4}$ by \cite{ps}. 
It follows that 
\begin{align*}
\|\,|S|-V^*U|S|^{1/2}|T|^{1/2}\|_1
&\le \|\,|S|-V^*S\|_1+\|V^*U|S|^{1/2}\|_2
 \|\,|S|^{1/2}-|T|^{1/2}\|_2\\
&\le \|S-T\|_1+\|\, |S|-|T| \,\|_1 +2\e^{1/4}\\
&\le \e+2\e^{1/2}+2\e^{1/4}.
\end{align*}
Therefore, we have 
\[
\|\p(S)-\p(T)\|_2^2=2-2\Re\Tr(\p(T)^*\p(S))
\le 2\e+4\e^{1/2}+4\e^{1/4}.
\]
This proves the uniform continuity of $\p$. 
The uniform continuity of $\p^{-1}$ follows from 
the Cauchy-Schwarz inequality. 
\end{proof}

Let us identify $S_1$ (resp.\ $S_2$) with 
the projective (resp.\ $\ell_2$-) tensor product 
$\ell_2\hat{\otimes}\ell_2$ (resp.\ $\ell_2\otimes_2\ell_2$).
The identification is given by 
$\xi\otimes\eta\mapsto T_{\xi,\eta}$, 
where 
$T_{\xi,\eta}(\zeta)=(\zeta,\bar{\xi})\eta$. 
(The symbol $\bar{\xi}$ means the complex conjugate of $\xi\in\ell_2$.)
Then, the noncommutative Mazur map 
$\p\colon\ell_2\hat{\otimes}\ell_2\to\ell_2\otimes_2\ell_2$ 
satisfies $\p\circ(U\otimes\bar{U})=(U\otimes\bar{U})\circ\p$ 
for any unitary operator $U$ and is uniformly homeomorphic on the spheres.
\section{Proof of Non-amenability}
\begin{proof}
Let $A$ be one of the Banach algebras listed 
in Theorem \ref{main}. 
Set $q$ be so that $1/p+1/q=1$. 
To give a finite set $F\subset A$ which violates 
the definition of amenability, we fix a finite set $\Sg$ 
of generators of $\G$. 
For notational simplicity, we set $\Sg^+=\Sg\cup\{g_\upsilon\}$ 
and $\pi^p(g_\upsilon)=\upsilon^p$. 
Then, 
\[
F=\{\pi^p(g) : g\in\Sg^+\} 
\subset\prod\B(\ell_p(\Ll))\subset \B(\ell_p(\Lambda))
\]
(cf.\ Section \ref{sl}) is a finite subset in $A$. 
We choose $\e>0$ later. 
Suppose there exists 
$T=\sum_{i=1}^r a_i\otimes b_i\in A\otimes A$ 
which satisfy the conditions (i) and (ii) 
in Definition \ref{def}. 

We fix a prime number $l$ until the very last step. 
For each $m\in\Lambda$, 
we define $T_l(m)\in\ell_2(\Ll)\hat{\otimes}\ell_2(\Ll)$ 
by 
\[
T_l(m)=\sum_{i=1}^r 
P_la_i\delta_m\otimes P_l^*b_i^*\delta_m^*
\in\ell_2(\Ll)\hat{\otimes}\ell_2(\Ll),
\]
where $P_l$ is the (orthogonal) projection 
from $\ell_p(\Lambda)$ onto $\ell_p(\Ll)$. 
By Lemma \ref{lpineq} (and the following remark), we have 
$\sum_{m\in\Lambda}\| T_l(m)\|_{\ell_2(\Ll)\hat{\otimes}\ell_2(\Ll)}
\le |\Ll|\,\|T\|_{A\hat{\otimes}A}$. 
Similarly 
\begin{align*}
\sum_{m\in\Lambda}\|T_l(m) & - (\pi^2_l(g)\otimes\pi^2_l(g))T_l(m)
 \|_{\ell_2(\Ll)\hat{\otimes}\ell_2(\Ll)}\\
&\le |\Ll|\,\| T-\pi^p(g)\cdot T\cdot\pi^p(g)^{-1}
 \|_{A\hat{\otimes}A}
\le|\Ll|\,\e
\end{align*}
for every $g\in\Sg^+$. 
On the other hand, since $\sum_{i=1}^r a_ib_i=1$, we have 
\[
\sum_{m\in\Lambda}\| T_l(m) \|_{\ell_2(\Ll)\hat{\otimes}\ell_2(\Ll)}
\geq \sum_{m\in\Lambda}
 \langle P_la_i\delta_m,\, P_l^*b_i^*\delta_m^* \rangle
=\Tr(P_l)=|\Ll|. 
\]
Therefore, there exists $m_l\in\Lambda$ such that 
\[
\|T_l(m_l) - (\pi^2_l(g)\otimes\pi^2_l(g))T_l(m_l)
 \|_{\ell_2(\Ll)\hat{\otimes}\ell_2(\Ll)}
\le \e|\Sg^+|\,\| T_l(m_l)\|_{\ell_2(\Ll)\hat{\otimes}\ell_2(\Ll)}
\]
for all $g\in\Sg^+$. 
Set $S_l=T_l(m_l)/\|T_l(m_l)\|_{\wedge}$. 
The application of the noncommutative Mazur map $\p$ to this inequality 
yields (cf.\ the last remark in Section \ref{mazur}) that 
\[
\|\p(S_l) - (\pi^2_l(g)\otimes\pi^2_l(g))\p(S_l)
\|_{\ell_2(\Ll)\otimes_2\ell_2(\Ll)}\le\omega_\p(\e|\Sg^+|)
=:\delta_0 
\]
for all $g\in\Sg^+$, where $\omega_{\p}$ is 
the modulus of continuity of $\p$ (cf.\ Appendix). 
Now, since $\pi^2_l\otimes\pi^2_l$ 
on $\ell_2(\Ll)\otimes_2\ell_2(\Ll)$ 
is a unitary representation of $\G$, 
Theorem \ref{prt} is applicable. 
It follows that there exists 
a $(\pi^2_l\otimes\pi^2_l)$-invariant 
vector $X$ such that 
\[
\| \p(S_l) -X \|_{\ell_2(\Ll)\otimes_2\ell_2(\Ll)}
\le R\delta_0, 
\]
where $R=R(\G,\Sg)$ is the constant appearing in Theorem \ref{prt}. 
But since $\G$ action on $\Ll$ is 2-transitive, 
$X=\lambda I + \mu E$ for some $\lambda,\mu\in\C$, 
where $I=|\Ll|^{-1/2}\sum_{s\in\Ll}\delta_s\otimes\delta_s$ 
and $E=|\Ll|^{-1}\sum_{s,t\in\Ll}\delta_s\otimes\delta_t$. 
It follows that 
\begin{align*}
(2-2|\Ll|^{-2})^{1/2}|\mu| 
&= |\mu|\,\|E- (\upsilon^2\otimes\upsilon^2)E
 \|_{\ell_2(\Ll)\otimes_2\ell_2(\Ll)}\\
&= \|X- (\upsilon^2\otimes\upsilon^2)X
 \|_{\ell_2(\Ll)\otimes_2\ell_2(\Ll)}\\
&\le 2R\delta_0+\|\p(S_l)-(\upsilon^2\otimes\upsilon^2)\p(S_l)
 \|_{\ell_2(\Ll)\otimes_2\ell_2(\Ll)}\\
&\le (2R+1)\delta_0.
\end{align*}
Therefore, 
\[
\| \p(S_l)-\lambda I\|_{\ell_2(\Ll)\otimes_2\ell_2(\Ll)}
\le (3R+1)\delta_0=:\delta_1. 
\]
and in particular $|\lambda| \geq \|\p(S_l)\|-\delta_1=1-\delta_1$.
On the other hand, 
since $\p(S_l)$ is of rank~$\le r$ 
(as a Hilbert-Schmidt operator on $\ell_2(\Ll)$), 
we have 
\[
(1-\frac{r}{|\Ll|})|\lambda|^2 
\le \|\p(S_l)-\lambda I\|_{\ell_2(\Ll)\otimes_2\ell_2(\Ll)}^2=\delta_1^2.
\]
Consequently, 
\[
r \geq (1-\frac{\delta_1^2}{(1-\delta_1)^2})|\Ll|
=:(1-\delta(\e))|\Ll|, 
\]
where $\delta\colon\R_+\to[0,1]$ is a function which does not 
depend on $l$ (but depends on $\G,\Sg$ and $\omega_\p$) 
such that $\lim_{\e\to 0}\delta(\e)=0$. 
Since $l$ can be arbitrarily large, 
this is absurd when $\e>0$ is chosen sufficiently small 
so that $\delta(\e)<1$.
\end{proof}
\begin{rem}
In the proof, one may use Corollary \ref{fga}, 
instead of Theorem~\ref{prt}, to deduce that 
$S_l\in\ell_2(\Ll)\hat{\otimes}\ell_2(\Ll)$ is 
close to a $(\pi_l^2\otimes\pi_l^2)$-invariant vector. 
The first half of the proof goes 
\textit{mutatis mutandis} for $B(\ell_p)$ and yields 
$S_l\in\ell_q(\Ll)\hat{\otimes}\ell_p(\Ll)$
which is almost $(\pi_l^q\otimes\pi_l^p)$-invariant. 
However, since we do not know whether the unit ball of 
$\ell_q\hat{\otimes}\ell_p$ is uniformly embeddable 
into a Hilbert space, it is not clear whether 
Corollary~\ref{fga} is applicable in this case or not. 
We note that it is not even known whether 
$\ell_q\hat{\otimes}\ell_p$ has non-trivial cotype or not 
(cf.\ \cite{pisier1}). 
\end{rem}
\appendix\section{Non-uniform Embeddability of Expanders}
In this appendix, we prove a result about 
embedding of expanders into a Banach space. 
Throughout this appendix, we denote by $G$ 
a graph which is simple, connected and $k$-regular. 
We abuse the notation and identify $G$ with its vertices. 
Further, we regard $G$ as the metric space equipped 
with the graph distance. 
For $A\subset G$, the boundary $\partial A$ is the set of 
(unoriented) edges which connects $A$ to $G\setminus A$. 
Let us recall some basic facts about expanders. 
A general reference for expanders is \cite{lubotzky}. 
The expanding constant $h(G)$ of $G$ is given by 
\[
h(G)=\inf\{ \frac{|\partial A|}{|A|} : A\subset V
\mbox{ with }|A|\leq|V|/2\}.
\]
For $h>0$, we say $G$ is an $(k,h)$-expander if 
it is $k$-regular and $h(G)\geq h$. 

Let $M_1$ and $M_2$ be metric spaces and let 
$\p\colon M_1\to M_2$ be a map. 
The Lipschitz norm of $\p$ is given by 
$\Lip(\p)=\sup\{\dist(\p(a),\p(b))/\dist(a,b) : a\neq b\}$. 
When $M_1=G$ is a graph, then 
\[
\Lip(\p)=\sup\{\dist(\p(s),\p(t)) : s,t\in G\mbox{ are adjacent}\}.
\]
The modulus of continuity of $\p$ is the function 
$\omega_\p\colon\R_+\to\R_+$ given by 
\[
\omega_\p(t)=\sup\{ \dist(\p(a),\p(b)) : 
a,b\in M_1\mbox{ and }\dist(a,b)\le t \}. 
\]
The map $\p$ is said to be uniformly continuous if 
$\lim_{t\to0}\omega_\p(t)=0$. 
The metric space $M_1$ is said to be 
uniformly homeomorphic to $M_2$ if 
there is a bijection $\p\colon M_1\to M_2$ 
such that both $\p$ and $\p^{-1}$ are uniformly continuous. 
Finally, $M_1$ is said to be 
uniformly embeddable into $M_2$ if 
$M_1$ is uniformly homeomorphic to a subset of $M_2$. 

We have the following concentration theorem 
for expanders embedded in a Banach space. 
A better estimate is known \cite{llr}\cite{matousek} 
for $X=\ell_p$ with $p\in[1,\infty)$. 

\begin{thm}\label{exp}
Let $X$ be a Banach space whose unit ball is 
uniformly embeddable into a Hilbert space. 
Then, for any $k\in\N$ and $h>0$, 
there exists a positive number $R=R(k,h,X)$ 
which satisfies the following; 
for any map $f$ from a $(k,h)$-expander $G$ 
into $X$, we have 
\[
\frac{1}{|G|}\sum_{s\in G}\|f(s)-m\|\le R\,\mathrm{Lip}(f),
\]
where $m=|G|^{-1}\sum_{s\in G}f(s)$ is the mean of $f$. 
\end{thm}

Any $n$-point metric space is isometric to 
a subset of $\ell_\infty^n$. 
A Banach space is said to have non-trivial cotype if 
it does not contain $\ell_\infty^n$'s uniformly. 
Thus a Banach space satisfying the conclusion of the theorem 
has non-trivial cotype (compare this with \cite{enflo}). 
It would be interesting to know whether 
the converse is also true or not. 
We note that, 
by Odell and Schlumprecht's theorem \cite{os}\cite{chaatit}, 
a Banach space $X$ with an unconditional basis 
and of non-trivial cotype has the unit ball 
which is uniformly homeomorphic to the unit ball of 
a Hilbert space. 
A general reference for related facts is \cite{bl}. 

We apply the above theorem to the Caley graphs of 
a finite group. 
Recall that the Caley graph $G(\Delta,\Sg)$ of a group $\Delta$ 
with a distinguished set of generators $\Sg$ is 
the graph whose vertices are elements of $\Delta$ 
and whose edges are pairs $\{g,h\}$ with $g^{-1}h\in\Sg$. 
The following is Margulis' celebrated construction 
of expanders (cf.\ \cite{lubotzky}). 
\begin{thm}\label{qcg}
Let $\G$ be a group with Kazhdan's property $\mathrm{(T)}$ 
and $\Sg$ be a finite set of generators of $\G$. 
Then, there exists $h=h(\G,\Sg)>0$ such that 
for any finite quotient group $\Delta$ of $\G$, 
the Caley graph $G(\Delta,\Sg)$ is an $(|\Sg|,h)$-expander. 
\end{thm}
\begin{cor}\label{fga}
Let $\G$ be a group with Kazhdan's property $\mathrm{(T)}$, 
$\Sg$ be a finite set of generators of $\G$ 
and let $X$ be a Banach space whose unit ball is 
uniformly embeddable into a Hilbert space. 
Then, there exists $R=R(\G,\Sg,X)>0$ 
which satisfies the following; 
for any representation $\pi\colon\G\to\B(X)$ 
with a finite image and any $\xi\in X$, 
there exists a $\pi$-invariant vector $\eta$ 
such that 
\[
\|\xi-\eta\|\le 
R\max_{g\in\G}\|\pi(g)\|^2\max_{h\in\Sg}\|\pi(h)\xi-\xi\|.
\]
\end{cor}
\begin{proof}
Apply Theorems \ref{exp} and \ref{qcg} 
to the map $f\colon\pi(\G)\ni\pi(g)\mapsto\pi(g)\xi\in X$.
\end{proof}

The following lemma is well-known, 
but we include the proof for completeness. 
\begin{lem}\label{con}
Let $G$ be a $(k,h)$-expander. Then, we have the following. 
\begin{enumerate}[(i)]
\item
For any map $f\colon G\to\ell_1$, we have 
\[
\frac{1}{|G|}\sum_{s\in G}\|f(s)-m\|
\leq 2\frac{k}{h}\Lip(f), 
\]
where $m=|G|^{-1}\sum_{s\in G}f(s)$ is the mean of $f$. 
\item
For any function $f\colon G\to\R_+$ with 
$|\{ s\in G : f(s)\leq R_0\}|\geq |G|/2$, 
we have
\[
\frac{1}{|G|}\sum_{s\in G}f(s)\leq R_0+\frac{k}{2h}\Lip(f).
\]
\end{enumerate}
\end{lem}
\begin{proof}
First, consider a function $g\colon G\to\R_+$ 
with $|\supp(g)|\le |G|/2$. 
Let $E$ be the set of unoriented edges of $G$ and 
let $A_r=\{ s\in G : g(s)\geq r\}$. 
We note that $|A_r| \le |\supp(g)| \le |G|/2$ and hence 
$|\partial A_r|\geq h|A_r|$ for $r>0$. 
Since $|g(s)-g(t)|$ is the Lebesgue measure of 
the interval between $g(s)$ and $g(t)$, 
we have 
\[
\sum_{\{s,t\}\in E}|g(s)-g(t)|
=\int_0^\infty|\partial A_r|\,\mathrm{d}r
\geq\int_0^\infty h|A_r|\,\mathrm{d}r
=h\sum_{s\in G}g(s). 
\]

Now the assertion (i) follows by considering each coordinate 
of $\ell_1$ and applying the above inequality to the positive 
(or perhaps negative) part of $\Re f$ (resp.\ $\Im f$). 
The part (ii) follows from the above inequality 
for $g=\max\{0,\,f-R_0\}$.
\end{proof}

Now, we give the proof of Theorem \ref{exp}. 

\begin{proof}[Proof of Theorem \ref{exp}]
We note that the Hilbert space $\ell_2$ 
is uniformly embeddable into $\ell_1$ (cf.\ \cite{bl}). 
Fix a uniform embedding $\p$ 
of the unit ball of $X$ into $\ell_1$. 
Take $R>0$ large enough so that it satisfies 
\[
R\geq 10\frac{k}{h}\quad\mbox{and}\quad
\omega_{\p^{-1}}(16\frac{k}{h}\omega_{\p}
(\frac{5}{R}))\le\frac{1}{9}.
\]
We proceed by contradiction. 
Suppose there exists a map 
$f\colon G\to X$ such that 
$(1/|G|)\sum_{s\in G}\|f(s)-m\|>R\,\Lip(f)$. 
We assume without loss of generality that 
$m=0$ and $\Lip(f)=1$. 
Choose $R_0>0$ such that 
there exists $s_0\in G$ satisfying
\[
|\{ x\in G : \|f(s)-f(s_0)\|\le R_0 \}|\geq\frac{3}{4}|G|.
\]
We may assume that $R_0$ attains its infimum under 
this condition. 
Since the mean of $f$ is zero, by Lemma \ref{con}, 
we have 
\[
R < \frac{1}{|G|}\sum_{s\in G}\|f(s)\| 
\le \frac{2}{|G|}\sum_{s\in G}\|f(s)-f(s_0)\|
\le 2R_0+\frac{k}{h} 
\]
and hence $R_0\geq (9/20)R\geq 9k/2h$. 
Let $\tilde{f}\colon G\to X$ be a map given by 
\[
\tilde{f}(x)=\begin{cases}
f(s) &\mbox{ if }\|f(s)-f(s_0)\|\le R_0\\ 
f(s_0)+R_0\frac{f(s)-f(s_0)}{\|f(s)-f(s_0)\|}
&\mbox{ if }\|f(s)-f(s_0)\|> R_0\end{cases}.
\] 
It follows that $\tilde{f}$ is 
a $2$-Lipschitz map from $G$ into the closed ball
of center $f(s_0)$ and radius $R_0$. 
Let $g\colon G\to \ell_1$ be the map 
given by $g(s)=\p((\tilde{f}(s)-f(s_0))/R_0)$. 
It follows $\Lip(g)\le\omega_{\p}(2/R_0)\le\omega_{\p}(5/R)$. 
Applying lemma \ref{con} to $g$, 
we obtain that 
\[
\frac{1}{|G|}\sum_{s\in G}\|g(s) - m(g)\|
\le 2\frac{k}{h}\omega_{\p}(\frac{5}{R})
=:\delta, 
\]
where $m(g)$ is the mean of $g$. 
It follows that for 
$A=\{ s\in G : \|g(s) - m(g) \|\le 4\delta \}$,
we have 
$|A|\geq(3/4)|G|$. 
Hence, for $A_1=A\cap\{ s\in G : f(s)=\tilde{f}(s)\}$, 
we have $|A_1|\geq|G|/2$. 
Choose $s_1\in A_1$. 
Note that for $s\in A_1$, we have 
$\| g(s) - g(s_1) \|\le 8\delta$ and hence 
$\| f(s) - f(s_1) \| \le \omega_{\p^{-1}}(8\delta)R_0\le R_0/9$ 
by the assumption on $R$. 
It follows from Lemma \ref{con} that 
\[
\frac{1}{|G|}\sum_{x\in G}\|f(s)-f(s_1)\| \le 
\frac{1}{9}R_0+\frac{k}{2h}\le \frac{2}{9}R_0.
\]
Therefore, we get 
\[
|\{ s\in G : \|f(s)-f(s_1)\| 
\le \frac{8}{9}R_0 \}|
\geq \frac{3}{4}|G|, 
\]
which is in contradiction with the minimality of $R_0$. 
\end{proof}


\begin{thebibliography}{LLR} 
%
\bibitem[Bh]{bhatia}
R. Bhatia, 
\textit{Matrix analysis}. 
Graduate Texts in Mathematics, \textbf{169}. 
Springer-Verlag, New York, 1997. 
%
\bibitem[BL]{bl}
Y. Benyamini and J. Lindenstrauss, 
\textit{Geometric nonlinear functional analysis}. Vol. 1. 
American Mathematical Society Colloquium Publications, \textbf{48}. 
American Mathematical Society, Providence, RI, 2000.
%
\bibitem[Ch]{chaatit}
F. Chaatit, 
\textit{On uniform homeomorphisms of the unit spheres 
of certain Banach lattices}. 
Pacific J. Math. \textbf{168} (1995), 11--31.
%
\bibitem[Co]{connes}
A. Connes, 
\textit{On the cohomology of operator algebras}. 
J. Funct. Anal. \textbf{28} (1978), 248--253. 
%
\bibitem[En]{enflo}
P. Enflo, 
\textit{On a problem of Smirnov}. 
Ark. Mat. \textbf{8} (1969), 107--109. 
%
\bibitem[HV]{hv}
P. de la Harpe and A. Valette, 
\textit{La propri{\'e}t{\'e} $(T)$ de Kazhdan pour les groupes 
localement compacts}, 
With an appendix by M. Burger. 
Ast{\'e}risque \textbf{175} (1989). 
%
\bibitem[Jo1]{johnson}
B. E. Johnson, 
\textit{Cohomology in Banach algebras}. 
Memoirs of the American Mathematical Society, No. \textbf{127}. 
American Mathematical Society, Providence, R.I., 1972.
%
\bibitem[Jo2]{johnson2}
B. E. Johnson, 
\textit{Approximate diagonals and cohomology of 
certain annihilator Banach algebras}. 
Amer. J. Math. \textbf{94} (1972), 685--698.
%
\bibitem[Ko]{kosaki}
H. Kosaki, 
\textit{On the continuity of the map $\varphi\to|\varphi|$ from 
the predual of a $W^*$-algebra}, 
J. Funct. Anal. \textbf{59} (1984), 123--131.
%
\bibitem[LLR]{llr}
N. Linial, E. London, Y. Rabinovich, 
\textit{The geometry of graphs and some of its algorithmic applications}. 
Combinatorica \textbf{15} (1995), 215--245.
%
\bibitem[Lu]{lubotzky}
A. Lubotzky, 
\textit{Discrete groups, expanding graphs and invariant measures}. 
With an appendix by J. D. Rogawski.

Progress in Mathematics, 125. 
Birkh{\"a}user Verlag, Basel, 1994.
%
\bibitem[Ma]{matousek}
J. Matou\v{s}ek, 
\textit{On embedding expanders into $\ell_p$ spaces}. 
Israel J. Math. \textbf{102} (1997), 189--197. 
%
\bibitem[OS]{os}
E. Odell and Th. Schlumprecht, 
\textit{The distortion problem}. 
Acta Math. \textbf{173} (1994), 259--281.
%
\bibitem[Pi1]{pisier1}
G. Pisier, 
\textit{Factorisation de fonctions analytiques \`{a} 
valeurs op\'{e}rateurs}. 
C. R. Acad. Sci. Paris S\'{e}r. I
Math. \textbf{307} (1988), 955--960.
%
\bibitem[Pi2]{pisier}
G. Pisier, 
\textit{On Read's proof that $B(\ell_1)$ is not amenable}. 
Seminar Report. 
%
\bibitem[PS]{ps}
R. Powers and E. St{\o}rmer, 
\textit{Free states of the canonical anticommutation relations}, 
Comm. Math. Phys. \textbf{16} (1970) 1--33.
%
\bibitem[Re]{read}
C. J. Read, 
\textit{Relative amenability and 
the non-amenability of $B(\ell^1)$}. 
Preprint. 
%
\bibitem[Ru]{runde}
V. Runde, 
\textit{Lectures on amenability}. 
Lecture Notes in Mathematics, \textbf{1774}. 
Springer-Verlag, Berlin, 2002.
%
\end{thebibliography}
\end{document}